\documentclass[11pt]{article}
\usepackage{amsfonts}
\usepackage{amsfonts}
\usepackage{amsfonts}
\usepackage{mathrsfs}
\usepackage{bm}
\usepackage{cite}
\usepackage[colorlinks,linkcolor=blue,citecolor=blue]{hyperref}
\usepackage{amssymb,amsmath} 
 \allowdisplaybreaks

\catcode`!=11
\let\!int\int \def\int{\displaystyle\!int}
\let\!lim\lim \def\lim{\displaystyle\!lim}
\let\!sum\sum \def\sum{\displaystyle\!sum}
\let\!sup\sup \def\sup{\displaystyle\!sup}
\let\!inf\inf \def\inf{\displaystyle\!inf}
\let\!cap\cap \def\cap{\displaystyle\!cap}
\let\!max\max \def\max{\displaystyle\!max}
\let\!min\min \def\min{\displaystyle\!min}
\let\!frac\frac \def\frac{\displaystyle\!frac}
\catcode`!=12

\let\oldsection\section
\renewcommand\section{\setcounter{equation}{0}\oldsection}

\allowdisplaybreaks
\def\pf{\it{Proof.}\rm\quad}

\def\N{\mathbb{N}}

\newtheorem{thm}{Theorem}[section]
\newtheorem{lem}[thm]{Lemma}
\newtheorem{cor}[thm]{Corollary}

\begin{document}
\title {\bf  Euler sums of generalized hyperharmonic numbers}
\author{
{Ce Xu\thanks{Corresponding author  Email: 15959259051@163.com (C. Xu)}}\\[1mm]
\small  School of Mathematical Sciences, Xiamen University\\
\small Xiamen
361005, P.R. China}
\date{}
\maketitle \noindent{\bf Abstract }  The generalized hyperharmonic numbers $h_n^{(m)}(k)$ are defined by means of
the multiple harmonic numbers. We show that the hyperharmonic numbers $h_n^{(m)}(k)$ satisfy certain recurrence relation which allow us to write them in terms of classical harmonic numbers. Moreover, we prove that the Euler-type
sums with hyperharmonic numbers:
\[S\left( {k,m;p} \right): = \sum\limits_{n = 1}^\infty  {\frac{{h_n^{\left( m \right)}\left( k \right)}}{{{n^p}}}} \;\;\left(p\geq m+1,\ {k = 1,2,3} \right)\]
can be expressed as a rational linear combination of products of Riemann zeta values and harmonic numbers. This is an extension of the results of Dil (2015) \cite{AD2015} and Mez$\ddot{\rm o}$ (2010) \cite{M2010}.  Some interesting new consequences and illustrative examples are considered.
\\[2mm]
\\[2mm]
\section{Introduction}
Let $\N:=\{1,2,3,\ldots\}$ be the set of natural numbers, $\N_0:=\N\cup \{0\}$, and $\mathbb{N} \setminus \{1\}:=\{2,3,4,\cdots\}$. Hyperharmonic numbers and their generalizations are classically defined by
\begin{align*}
&h_n^{\left( m \right)}\left( k \right): = \sum\limits_{\scriptstyle 1 \le {n_{m + k - 1}} <  \cdots  < {n_m} \hfill \atop
  \scriptstyle  \le {n_{m - 1}} \le  \cdots  \le {n_{1}} \le n \hfill} {\frac{1}{{{n_m}{n_{m + 1}} \cdots {n_{m + k - 1}}}}}, \tag{1.1}\\
&h_n^{\left( m \right)}\left( 1 \right) \equiv h_n^{\left( m \right)}: = \sum\limits_{1 \le {n_m} \le  \cdots  \le {n_{1}} \le n} {\frac{1}{{{n_m}}}} ,\tag{1.2}
\end{align*}
where $k,m,n\in \N$ and for any $n<k$, we set $h_n^{\left( m \right)}\left( k \right): = 0$. When $k=1$ in (1.1), the number $h_n^{\left( m \right)}\left( 1 \right) \equiv h_n^{\left( m \right)}$ is called the classical hyperharmonic number (see \cite{CG1996,AD2009,AD2008,AD2015,M2010}).  In special, the hyperharmonic number $h_n^{\left( 1 \right)}$ is
simply called the classical harmonic number, which is the sum of the reciprocals of the first $n$ natural numbers:
\[h_n^{\left( 1 \right)} \equiv {H_n}: = \sum\limits_{k = 1}^n {\frac{1}{k}} .\]
Moreover, in \cite{M2010}, Mez$\ddot{\rm o}$ and Dil showed that $h_n^{\left( m \right)}$ can be expressed by binomial coefficients and classical harmonic numbers:
\[h_n^{\left( m \right)} = \left( {\begin{array}{*{20}{c}}
   {n + m - 1}  \\
   {m - 1}  \\
\end{array}} \right)\left( {{H_{n + m - 1}} - {H_{m - 1}}} \right).\]
The $n$-th generalized harmonic numbers of order $k$, denoted by $H^{(k)}_n$, is defined by
\[\ H^{(k)}_n:=\sum\limits_{j=1}^n\frac {1}{j^k},\ n, k \in \N,\tag{1.3}\] where the empty sum $H^{(k)}_0$ is conventionally understood to be zero, and $H^{(1)}_n\equiv H_n$. The limit as $n$ tends to infinity exists if $k > 1$.
In the limit of $n\rightarrow \infty$, the generalized harmonic number converges to the Riemann zeta value:
 \[\mathop {\lim }\limits_{n \to \infty } H_n^{\left( k \right)} = \zeta( k),\ {\Re} \left( k \right) > 1,\ k\in \N,\]
where the Riemann zeta function is defined by
\[\zeta(s):=\sum\limits_{n = 1}^\infty {\frac {1}{n^{s}}},\Re(s)>1.\tag{1.4}\]
In general, for $r\in \N$, ${\bf s}:=(s_1,s_2,\ldots,s_r)\in \N^r$, and a non-negative integer $n$, the multiple harmonic number is defined by
\[H_n^{\left( {{s_1},{s_2}, \cdots ,{s_r}} \right)}: = \sum\limits_{1 \le {n_r} < {n_{r - 1}} <  \cdots  < {n_1} \le n} {\frac{1}{{{n_1^{{s_1}}}{n_2^{{s_2}}} \cdots {n_r^{{s_r}}}}}}.\tag{1.5}\]
By convention, we put $H^{(\bf s)}_n=0$, if $n<r$, and $H_n^{\left( \emptyset  \right)} = 1$. The limit cases of multiple harmonic numbers give rise to multiple zeta values:
\[\zeta \left( {{s_1},{s_2}, \cdots ,{s_r}} \right) = \mathop {\lim }\limits_{n \to \infty } H_n^{\left( {{s_1},{s_2}, \cdots ,{s_r}} \right)}\]
defined for $s_2,s_3,\ldots,s_r\geq1$ and $s_1\geq2$ to ensure convergence of the series. Here, $w:={s_1} +  \cdots  + {s_r}$ and $r$ are called the weight and the multiplicity, respectively.
To simplify the reading of such formulas, when a string of arguments is repeated an exponent is used. In other words, we treat string multiplication as concatenation. For example,
\[H_n^{\left( {\underbrace {1, \cdots ,1}_r} \right)} = H_n^{\left( {{{\left\{ 1 \right\}}^r}} \right)},\;H_n^{\left( {\underbrace {2, \cdots ,2}_p,\underbrace {3, \cdots ,3}_r} \right)} = H_n^{\left( {{{\left\{ 2 \right\}}^p},{{\left\{ 3 \right\}}^r}} \right)}.\]
With these notations, then the definition of hyperharmonic number $h_n^{\left( m \right)}\left( k \right)$ of formula (1.1) can be rewritten as
\[h_n^{\left( m \right)}\left( k \right): = \sum\limits_{1 \le {n_m} \le {n_{m - 1}} \le  \cdots  \le {n_{1}} \le n} {\frac{{H_{{n_m} - 1}^{\left( {{{\left\{ 1 \right\}}^{k - 1}}} \right)}}}{{{n_m}}}},\tag{1.6} \]
where ${{H_{{n_m} - 1}^{\left( {{{\left\{ 1 \right\}}^{k }}} \right)}}}$ is the multiple harmonic number  ${{H_{n}^{\left( {{{\left\{ 1 \right\}}^{k}}} \right)}}}$ with $n=n_m-1$.

The subject of this paper is Euler-type sums $S\left( {k,m;p} \right)$, which is the infinite sum whose general term is a product of hyperharmonic numbers and a power of $n^{-1}$. Here, $p> m$ is both necessary and sufficient for the sum $S\left( {k,m;p} \right)$ to converge. The classical linear Euler sum is defined by (\cite{FS1998})
\[{S_{p,q}} := \sum\limits_{n = 1}^\infty  {\frac{{{H^{(p)} _n}}}{{{n^q}}}},\ \ p\in \N, q\in \mathbb{N} \setminus \{1\}.\tag{1.7} \]
The number $w=p+q$ is defined as the weight of $S_{p,q}$. The evaluation of $S_{p,q}$ in terms of values of Riemann zeta function at positive integers is known when $p=1,\ p=q,\ (p,q)=(2,4),(4,2)$ or $p+q$ is odd (see \cite{BBG1994,BBG1995,BZB2008,FS1998}). For example, Euler discovered the following formula
\[S_{1,k}=\sum\limits_{n = 1}^\infty  {\frac{{{H_n}}}{{{n^k}}}}  = \frac{1}{2}\left\{ {\left( {k + 2} \right)\zeta \left( {k + 1} \right) - \sum\limits_{i = 1}^{k - 2} {\zeta \left( {k - i} \right)\zeta \left( {i + 1} \right)} } \right\}.\tag{1.8}\]
Related series were studied by Mez$\ddot{\rm o}$ in \cite{M2014}, Sofo \cite{S2015} and Xu.et al \cite{Xu2016,X2016}, for instance. Similarly, it has been discovered in the course of the years that many Euler type sums $S\left( {k,m;p} \right)$ admit expressions involving finitely the zeta values, that is to say values of the Riemann zeta function at the positive integer arguments, for more details, see for instance \cite{AD2015,M2010}. For example, Dil and Boyadzhiev \cite{AD2015} gave explicit reductions to zeta values and (unsigned) Stirling numbers of the first kind for all sums $S\left( {k,m;p} \right)$ with $k=1$. Here, the (unsigned) Stirling number of the first kind $\left[ {\begin{array}{*{20}{c}}
   n  \\
   k  \\
\end{array}} \right]$ is defined by \cite{CG1996,L1974}
\[n!x\left( {1 + x} \right)\left( {1 + \frac{x}{2}} \right) \cdots \left( {1 + \frac{x}{n}} \right) = \sum\limits_{k = 0}^n {\left[ {\begin{array}{*{20}{c}}
   {n + 1}  \\
   {k + 1}  \\
\end{array}} \right]{x^{k + 1}}} \tag{1.9} \]
with $\left[ {\begin{array}{*{20}{c}}
   n  \\
   k  \\
\end{array}} \right]=0$, if $n<k$ and $\left[ {\begin{array}{*{20}{c}}
   n  \\
   0  \\
\end{array}} \right]=\left[ {\begin{array}{*{20}{c}}
   0  \\
   k  \\
\end{array}} \right]=0,\ \left[ {\begin{array}{*{20}{c}}
   0  \\
   0  \\
\end{array}} \right]=1$, or equivalently, by the generating function:
 \[{\log ^k}\left( {1 - x} \right) = {\left( { - 1} \right)^k}k!\sum\limits_{n = 1}^\infty  {\left[ {\begin{array}{*{20}{c}}
   n  \\
   k  \\
\end{array}} \right]\frac{{{x^n}}}{{n!}}} ,\;x \in \left[ { - 1,1} \right).\tag{1.10} \]
Moreover, the (unsigned) Stirling numbers $\left[ {\begin{array}{*{20}{c}}
   n  \\
   k  \\
\end{array}} \right]$ of the first kind satisfy a recurrence relation in the form
\[\left[ {\begin{array}{*{20}{c}}
   n  \\
   k  \\
\end{array}} \right] = \left[ {\begin{array}{*{20}{c}}
   {n - 1}  \\
   {k - 1}  \\
\end{array}} \right] + \left( {n - 1} \right)\left[ {\begin{array}{*{20}{c}}
   {n - 1}  \\
   k  \\
\end{array}} \right].\tag{1.11} \]
By the definition of $\left[ {\begin{array}{*{20}{c}}
   n  \\
   k  \\
\end{array}} \right]$, we see that we may rewrite (1.9) as
\begin{align*}
\sum\limits_{k = 0}^n {\left[ {\begin{array}{*{20}{c}}
   {n + 1}  \\
   {k + 1}  \\
\end{array}} \right]}x^{k}&= n!\exp \left\{ {\sum\limits_{j = 1}^n {\ln \left( {1 + \frac{x}{j}} \right)} } \right\}\\
& = n!\exp \left\{ {\sum\limits_{j = 1}^n {\sum\limits_{k = 1}^\infty  {{{\left( { - 1} \right)}^{k - 1}}\frac{{{x^k}}}{{k{j^k}}}} } } \right\}\\
& = n!\exp \left\{ {\sum\limits_{k = 1}^\infty  {{{\left( { - 1} \right)}^{k - 1}}\frac{{{H^{(k)} _n}{x^k}}}{k}} } \right\}.
\end{align*}
Therefore, we know that ${\left[ {\begin{array}{*{20}{c}}
   n  \\
   k  \\
\end{array}} \right]}$ is a rational linear combination of products of harmonic numbers. Moreover, we deduce the following identities
\begin{align*}
&\left[ {\begin{array}{*{20}{c}}
   n  \\
   1  \\
\end{array}} \right]= \left( {n - 1} \right)!,\left[ {\begin{array}{*{20}{c}}
   n  \\
   2  \\
\end{array}} \right] = \left( {n - 1} \right)!{H_{n - 1}},\left[ {\begin{array}{*{20}{c}}
   n  \\
   3 \\
\end{array}} \right] = \frac{{\left( {n - 1} \right)!}}{2}\left[ {H_{n - 1}^2 - {H^{(2)} _{n - 1}}} \right],\\
&\left[ {\begin{array}{*{20}{c}}
   n  \\
   4  \\
\end{array}} \right] = \frac{{\left( {n - 1} \right)!}}{6}\left[ {H_{n - 1}^3 - 3{H_{n - 1}}{H^{(2)} _{n - 1}} + 2{H^{(3)} _{n - 1}}} \right], \\
&\left[ {\begin{array}{*{20}{c}}
   n  \\
   5  \\
\end{array}} \right] = \frac{{\left( {n - 1} \right)!}}{{24}}\left[ {H_{n - 1}^4 - 6{H^{(4)} _{n - 1}} - 6H_{n - 1}^2{H^{(2)} _{n - 1}} + 3(H^{(2)} _{n - 1})^2 + 8H_{n - 1}^{}{H^{(3)} _{n - 1}}} \right].
\end{align*}
In this paper we are interested in Euler-type sums with hyperharmonic numbers $S\left( {k,m;p} \right)$. Such series could be of interest in analytic number theory. We will prove that the generalized hyperharmonic number $h_n^{\left( m \right)}\left( k \right)$ can be expressed as a rational linear combination of products of harmonic numbers. Furthermore, we also provide an explicit evaluation of $S\left( {k,m;p} \right)$ with $k=2,3$ in a closed form in terms of zeta values and Stirling numbers
of the first kind. The results which we present here can be seen as an extension of Mez$\ddot{\rm o}$ and Dil's work.

\section{Main Theorems and their Proof}

In this section, we will show that the hyperharmonic number $h_n^{\left( m \right)}\left( k \right)$ is expressible in terms of harmonic numbers and give recurrence formula. We need the following lemma.
\begin{lem} For positive integers $n$ and $k$, then the following identity holds:
\[\left[ {\begin{array}{*{20}{c}}
   n  \\
   k  \\
\end{array}} \right] = \left( {n - 1} \right)!H_{n - 1}^{\left( {{{\left\{ 1 \right\}}^{k - 1}}} \right)}.\tag{2.1}\]
\end{lem}
\pf By considering the generating function (1.10), we know that we need to prove the following identity
\[{\log ^k}\left( {1 - x} \right) = {\left( { - 1} \right)^k}k!\sum\limits_{n = 1}^\infty  {H_{n - 1}^{\left( {{{\left\{ 1 \right\}}^{k - 1}}} \right)}\frac{{{x^n}}}{n}}.\tag{2.2}\]
To prove the identity we proceed by induction on $k$. Obviously, it is valid for $k=1$. For $k>1$ we use the integral identity
\[{\log ^{k{\rm{ + }}1}}\left( {1 - x} \right){\rm{ = }} - \left( {k + 1} \right)\int\limits_0^x {\frac{{{{\log }^k}\left( {1 - t} \right)}}{{1 - t}}dt} \]
and apply the induction hypothesis, by using Cauchy product of power series, we arrive at
\begin{align*}
{\log ^{k{\rm{ + }}1}}\left( {1 - x} \right){\rm{ = }}& - \left( {k + 1} \right)\int\limits_0^x {\frac{{{{\log }^k}\left( {1 - t} \right)}}{{1 - t}}dt} \\
& = {\left( { - 1} \right)^{k + 1}}\left( {k + 1} \right)!\sum\limits_{n = 1}^\infty  {\frac{1}{{n + 1}}\sum\limits_{i = 1}^n \frac {{H_{i - 1}^{\left( {{{\left\{ 1 \right\}}^{k - 1}}} \right)}}}{i} } {x^{n + 1}}\\
& = {\left( { - 1} \right)^{k + 1}}\left( {k + 1} \right)!\sum\limits_{n = 1}^\infty  {\frac{{H_{n}^{\left( {{{\left\{ 1 \right\}}^{k}}} \right)}}{}}{{n + 1}}} {x^{n + 1}}.
\end{align*}
Nothing that $ {{H_{n}^{\left( {{{\left\{ 1 \right\}}^{k}}} \right)}}{}}= 0$ when $n<k$. Hence,
we can deduce (2.2) holds. Thus, comparing the coefficients of $x^n$ in (1.10) and (2.2), we obtain formula (2.1). The proof of lemma 2.1 is completed.\hfill$\square$\\
By using (1.6) and (2.2), we find that the generating function of hyperharmonic number $h_n^{\left( m \right)}\left( k \right)$ is given as
\[\sum\limits_{n = 1}^\infty  {h_n^{\left( m \right)}\left( k \right){z^n}}  = \frac{{{{\left( { - 1} \right)}^k}}}{{k!}}\frac{{{{\log }^k}\left( {1 - z} \right)}}{{{{\left( {1 - z} \right)}^m}}},\;z \in \left[ { - 1,1} \right).\tag{2.3}\]
On the other hand, we note that the function on the right hand side of (2.3) is equal to
\[\frac{{{{\left( { - 1} \right)}^k}}}{{k!}}\frac{{{{\log }^k}\left( {1 - z} \right)}}{{{{\left( {1 - z} \right)}^m}}} = \frac{1}{{k!}}\mathop {\lim }\limits_{x \to m} \frac{{{\partial ^k}}}{{\partial {x^k}}}\left( {\frac{1}{{{{\left( {1 - z} \right)}^x}}}} \right)\ (k,m\in \N_0).\tag{2.4}\]
Therefore, the relations (2.3) and (2.4) yield the following result:
\[\sum\limits_{n = 1}^\infty  {h_n^{\left( m \right)}\left( k \right){z^n}}  = \frac{1}{{k!}}\mathop {\lim }\limits_{x \to m} \frac{{{\partial ^k}}}{{\partial {x^k}}}\left( {\frac{1}{{{{\left( {1 - z} \right)}^x}}}} \right).\tag{2.5}\]
Moreover, we know that the generating function of ${\left( {1 - z} \right)^{ - x}}$ is given as
\[\frac{1}{{{{\left( {1 - z} \right)}^x}}} = \sum\limits_{n = 0}^\infty  {\frac{{{{\left( x \right)}_n}}}{{n!}}{z^n}} ,\;z \in \left( { - 1,1} \right),\tag{2.6}\]
where ${{{\left( x \right)}_n}}$ represents the Pochhammer symbol (or the shifted factorial) given by
\[{\left( x \right)_n} := x\left( {x + 1} \right) \cdots \left( {x + n - 1} \right)\tag{2.7}\]
with $(x)_0:=1$. Hence, upon differentiating both members of (2.6) $k$ times with respect to $x$ then setting $x=m$, and combining (2.5), we readily arrive at the following relationship:
\[h_n^{\left( m \right)}\left( k \right) = \frac{1}{{k!n!}}\mathop {\lim }\limits_{x \to m} \frac{{{\partial ^k}{{\left( x \right)}_n}}}{{\partial {x^k}}},\ k\in \N\tag{2.8}.\]
By convention, from (2.8), we define that
\[h_n^{\left( m \right)}\left( 0 \right): = \frac{1}{{n!}}{\left( m \right)_n} = \left( {\begin{array}{*{20}{c}}
   {m + n - 1}  \\
   {m - 1}  \\
\end{array}} \right).\tag{2.9}\]
By simple calculation, the $\frac{{{\partial ^k}{{\left( x \right)}_n}}}{{\partial {x^k}}}$ satisfy a recurrence relation in the form
\[\frac{{{\partial ^k}{{\left( x \right)}_n}}}{{\partial {x^k}}} = \sum\limits_{i = 0}^{k - 1} {\left( {\begin{array}{*{20}{c}}
   {k - 1}  \\
   i  \\
\end{array}} \right)\frac{{{\partial ^i}{{\left( x \right)}_n}}}{{\partial {x^i}}}} \left[ {{\psi ^{\left( {m -i - 1} \right)}}\left( {x + n} \right) - {\psi ^{\left( {m - i- 1} \right)}}\left( x \right)} \right],\ k\in \N.\tag{2.10}\]
Here, ${\psi ^{\left( m \right)}}\left( x \right)$ stands for the polygamma function of order $m$ defined as the $(m+1)$ th derivative of the logarithm of the gamma function:
\[{\psi ^{\left( m \right)}}\left( x \right): = \frac{{{d^m}}}{{d{x^m}}}\psi \left( x \right) = \frac{{{d^{m+1}}}}{{d{x^{m+1}}}}\log \Gamma \left( x \right).\]
Thus 
\[{\psi ^{\left( 0 \right)}}\left( x \right) = \psi \left( x \right) = \frac{{\Gamma '\left( x \right)}}{{\Gamma \left( x \right)}}\]
holds where $\psi (x)$ is the digamma function and $\Gamma \left( x \right)$ is the gamma function.
${\psi ^{\left( m \right)}}\left( x \right)$ satisfy the following relations
\[\psi \left( z \right) =  - \gamma  + \sum\limits_{n = 0}^\infty  {\left( {\frac{1}{{n + 1}} - \frac{1}{{n + z}}} \right)} ,\;z\notin  \N^-_0:=\{0,-1,-2\ldots\}, \]
\[{\psi ^{\left( n \right)}}\left( z \right) = {\left( { - 1} \right)^{n + 1}}n!\sum\limits_{k = 0}^\infty  {1/{{\left( {z + k} \right)}^{n + 1}}}, n\in \N,\]
\[\psi \left( {x + n} \right) = \frac{1}{x} + \frac{1}{{x + 1}} +  \cdots  + \frac{1}{{x + n - 1}} + \psi \left( x \right),\;n \in \N .\]
Here, $\gamma$ denotes the Euler-Mascheroni constant, defined by
\[\gamma  := \mathop {\lim }\limits_{n \to \infty } \left( {\sum\limits_{k = 1}^n {\frac{1}{k}}  - \ln n} \right) =  - \psi \left( 1 \right) \approx {\rm{ 0 }}{\rm{. 577215664901532860606512 }}....\]
Hence, combining (2.8), (2.9) and (2.10), we obtain the recurrence relation
\[h_n^{\left( m \right)}\left( k \right) = \frac{{{{\left( { - 1} \right)}^{k - 1}}}}{k}\sum\limits_{i = 0}^{k - 1} {{{\left( { - 1} \right)}^i}h_n^{\left( m \right)}\left( i \right)\left\{ {H_{m + n - 1}^{\left( {k - i} \right)} - H_{m - 1}^{\left( {k - i} \right)}} \right\}} .\tag{2.11}\]
By (2.11), we give the following description of hyperharmonic number $h_n^{\left( m \right)}\left( k \right)$.
\begin{thm} For positive integers $n$ and $k$, then the hyperharmonic number $h_n^{\left( m \right)}\left( k \right)$ can be expressed in terms of ordinary harmonic numbers.
\end{thm}
For example, setting $m =1, 2, 3, 4$ in the above equation (2.11) we obtain
\begin{align*}
&h_n^{\left( m \right)}\left( 1 \right) = \left( {\begin{array}{*{20}{c}}
   {n + m - 1}  \\
   {m - 1}  \\
\end{array}} \right)\left( {{H_{n + m - 1}} - {H_{m - 1}}} \right),\tag{2.12}\\
&h_n^{\left( m \right)}\left( 2 \right) = \frac{1}{2}\left( {\begin{array}{*{20}{c}}
   {n + m - 1}  \\
   {m - 1}  \\
\end{array}} \right)\left\{ {{{\left( {{H_{n + m - 1}} - {H_{m - 1}}} \right)}^2} - \left( {H_{n + m - 1}^{\left( 2 \right)} - H_{m - 1}^{\left( 2 \right)}} \right)} \right\},\tag{2.13}\\
&h_n^{\left( m \right)}\left( 3 \right) = \frac{1}{{3!}}\left( {\begin{array}{*{20}{c}}
   {n + m - 1}  \\
   {m - 1}  \\
\end{array}} \right)\left\{ \begin{array}{l}
 {\left( {{H_{n + m - 1}} - {H_{m - 1}}} \right)^3} + 2\left( {H_{n + m - 1}^{\left( 3 \right)} - H_{m - 1}^{\left( 3 \right)}} \right) \\
  - 3\left( {{H_{n + m - 1}} - {H_{m - 1}}} \right)\left( {H_{n + m - 1}^{\left( 2 \right)} - H_{m - 1}^{\left( 2 \right)}} \right) \\
 \end{array} \right\},\tag{2.14}\\
 &h_n^{\left( m \right)}\left( 4 \right) = \frac{1}{{4!}}\left( {\begin{array}{*{20}{c}}
   {n + m - 1}  \\
   {m - 1}  \\
\end{array}} \right)\left\{ \begin{array}{l}
 {\left( {{H_{n + m - 1}} - {H_{m - 1}}} \right)^4} \\
  + 3{\left( {H_{n + m - 1}^{\left( 2 \right)} - H_{m - 1}^{\left( 2 \right)}} \right)^2} - 6\left( {H_{n + m - 1}^{\left( 4 \right)} - H_{m - 1}^{\left( 4 \right)}} \right) \\
  - 6{\left( {{H_{n + m - 1}} - {H_{m - 1}}} \right)^2}\left( {H_{n + m - 1}^{\left( 2 \right)} - H_{m - 1}^{\left( 2 \right)}} \right) \\
  + 8\left( {{H_{n + m - 1}} - {H_{m - 1}}} \right)\left( {H_{n + m - 1}^{\left( 3 \right)} - H_{m - 1}^{\left( 3 \right)}} \right) \\
 \end{array} \right\}.\tag{2.15}
\end{align*}
By replacing $x$ by $n$ and $n$ by $r$ in (1.9), we deduce that
\[\left( {\begin{array}{*{20}{c}}
   {n + r}  \\
   r  \\
\end{array}} \right) = \frac{1}{{r!}}\sum\limits_{k = 1}^{r + 1} {\left[ {\begin{array}{*{20}{c}}
   {r + 1}  \\
   k  \\
\end{array}} \right]{n^{k - 1}}} .\tag{2.16}\]
Therefore, the relations (2.13), (2.14) and (2.16) yield the following results
\begin{align*}
h_n^{\left( {r + 1} \right)}\left( 2 \right) =& \frac{1}{2}\left( {\begin{array}{*{20}{c}}
   {n + r}  \\
   r  \\
\end{array}} \right)\left\{ {{{\left( {{H_{n + r}} - {H_r}} \right)}^2} - \left( {H_{n + r}^{\left( 2 \right)} - H_r^{\left( 2 \right)}} \right)} \right\}\\
 =& \frac{1}{{2!r!}}\sum\limits_{k = 1}^{r + 1} {\left[ {\begin{array}{*{20}{c}}
   {r + 1}  \\
   k  \\
\end{array}} \right]{n^{k - 1}}\left\{ {{{\left( {{H_{n + r}} - {H_r}} \right)}^2} - \left( {H_{n + r}^{\left( 2 \right)} - H_r^{\left( 2 \right)}} \right)} \right\}},\tag{2.17}\\
h_n^{\left( {r + 1} \right)}\left( 3 \right) =& \frac{1}{{3!}}\left( {\begin{array}{*{20}{c}}
   {n + r}  \\
   r  \\
\end{array}} \right)\left\{ \begin{array}{l}
 {\left( {{H_{n + r}} - {H_r}} \right)^3} + 2\left( {H_{n + r}^{\left( 3 \right)} - H_r^{\left( 3 \right)}} \right) \\
  - 3\left( {{H_{n + r}} - {H_r}} \right)\left( {H_{n + r}^{\left( 2 \right)} - H_r^{\left( 2 \right)}} \right) \\
 \end{array} \right\}\\
  =& \frac{1}{{3!r!}}\sum\limits_{k = 1}^{r + 1} {\left[ {\begin{array}{*{20}{c}}
   {r + 1}  \\
   k  \\
\end{array}} \right]{n^{k - 1}}\left\{ \begin{array}{l}
 {\left( {{H_{n + r}} - {H_r}} \right)^3} + 2\left( {H_{n + r}^{\left( 3 \right)} - H_r^{\left( 3 \right)}} \right) \\
  - 3\left( {{H_{n + r}} - {H_r}} \right)\left( {H_{n + r}^{\left( 2 \right)} - H_r^{\left( 2 \right)}} \right) \\
 \end{array} \right\}} ,\tag{2.18}\\
h_n^{\left( {r + 1} \right)}\left( 4 \right) =& \frac{1}{{4!}}\left( {\begin{array}{*{20}{c}}
   {n + r}  \\
   r  \\
\end{array}} \right)\left\{ \begin{array}{l}
 {\left( {{H_{n + r}} - {H_r}} \right)^4} + 3{\left( {H_{n + r}^{\left( 2 \right)} - H_r^{\left( 2 \right)}} \right)^2} \\
  - 6\left( {H_{n + r}^{\left( 4 \right)} - H_r^{\left( 4 \right)}} \right) \\
  - 6{\left( {{H_{n + r}} - {H_r}} \right)^2}\left( {H_{n + r}^{\left( 2 \right)} - H_r^{\left( 2 \right)}} \right) \\
  + 8\left( {{H_{n + r}} - {H_r}} \right)\left( {H_{n + r}^{\left( 3 \right)} - H_r^{\left( 3 \right)}} \right) \\
 \end{array} \right\}\\
  = &\frac{1}{{4!r!}}\sum\limits_{k = 1}^{r + 1} {\left[ {\begin{array}{*{20}{c}}
   {r + 1}  \\
   k  \\
\end{array}} \right]{n^{k - 1}}\left\{ \begin{array}{l}
 {\left( {{H_{n + r}} - {H_r}} \right)^4} + 3{\left( {H_{n + r}^{\left( 2 \right)} - H_r^{\left( 2 \right)}} \right)^2} \\
  - 6\left( {H_{n + r}^{\left( 4 \right)} - H_r^{\left( 4 \right)}} \right) \\
  - 6{\left( {{H_{n + r}} - {H_r}} \right)^2}\left( {H_{n + r}^{\left( 2 \right)} - H_r^{\left( 2 \right)}} \right) \\
  + 8\left( {{H_{n + r}} - {H_r}} \right)\left( {H_{n + r}^{\left( 3 \right)} - H_r^{\left( 3 \right)}} \right) \\
 \end{array} \right\}}.\tag{2.19}
\end{align*}
Furthermore, using (2.17) and (2.18), by a direct calculation, we can give the following corollary.
\begin{cor} For integers $r\in \N_0$ and $n\in \N$, we have
\begin{align*}
&h_n^{\left( {r + 1} \right)}\left( 2 \right) = \frac{1}{{2!r!}}\sum\limits_{k = 1}^{r + 1} {\left[ {\begin{array}{*{20}{c}}
   {r + 1}  \\
   k  \\
\end{array}} \right]{n^{k - 1}}\left\{ {\left( {H_{n + r}^2 - H_{n + r}^{\left( 2 \right)}} \right) - 2{H_r}{H_{n + r}} + H_r^2 + H_r^{\left( 2 \right)}} \right\}},\tag{2.20} \\
&h_n^{\left( {r + 1} \right)}\left( 3 \right) = \frac{1}{{3!r!}}\sum\limits_{k = 1}^{r + 1} {\left[ {\begin{array}{*{20}{c}}
   {r + 1}  \\
   k  \\
\end{array}} \right]{n^{k - 1}}\left\{ \begin{array}{l}
 \left( {H_{n + r}^3 - 3{H_{n + r}}H_{n + r}^{\left( 2 \right)} + 2H_{n + r}^{\left( 3 \right)}} \right) - 3{H_r}\left( {H_{n + r}^2 - H_{n + r}^{\left( 2 \right)}} \right) \\
  + 3\left( {H_r^2 + H_r^{\left( 2 \right)}} \right){H_{n + r}} - \left( {H_r^3 + 3{H_r}H_r^{\left( 2 \right)} + 2H_r^{\left( 3 \right)}} \right) \\
 \end{array} \right\}}.\tag{2.21}
\end{align*}
\end{cor}
Moreover, from the definition of harmonic numbers $H^{(k)}_n$, we get
\[H_{n + r}^{\left( k \right)} = H_n^{\left( k \right)} + \sum\limits_{j = 1}^r {\frac{1}{{{{\left( {n + j} \right)}^k}}}},\ k,n\in \N. \tag{2.22}\]
By simple calculation, the following identities are easily derived
\begin{align*}
H_{n + r}^2 - H_{n + r}^{\left( 2 \right)} &= {\left( {{H_n} + \sum\limits_{j = 1}^r {\frac{1}{{n + j}}} } \right)^2} - \left( {H_n^{\left( 2 \right)} + \sum\limits_{j = 1}^r {\frac{1}{{{{\left( {n + j} \right)}^2}}}} } \right)\\
&=H_n^2 - H_n^{\left( 2 \right)} + 2{H_n}\left( {\sum\limits_{j = 1}^r {\frac{1}{{n + j}}} } \right) + 2\sum\limits_{1 \le i < j \le r} {\frac{1}{{\left( {n + i} \right)\left( {n + j} \right)}}},\tag{2.23}
\end{align*}
\begin{align*}
&H_{n + r}^3 - 3{H_{n + r}}H_{n + r}^{\left( 2 \right)} + 2H_{n + r}^{\left( 3 \right)} \\&= {\left( {{H_n} + \sum\limits_{j = 1}^r {\frac{1}{{n + j}}} } \right)^3} + 2\left( {H_n^{\left( 3 \right)} + \sum\limits_{j = 1}^r {\frac{1}{{{{\left( {n + j} \right)}^3}}}} } \right)\\
&\quad - 3\left( {{H_n} + \sum\limits_{j = 1}^r {\frac{1}{{n + j}}} } \right)\left( {H_n^{\left( 2 \right)} + \sum\limits_{j = 1}^r {\frac{1}{{{{\left( {n + j} \right)}^2}}}} } \right)\\
&=H_n^3 - 3{H_n}H_n^{\left( 2 \right)} + 2H_n^{\left( 3 \right)} + 3\left( {H_n^2 - H_n^{\left( 2 \right)}} \right)\left( {\sum\limits_{j = 1}^r {\frac{1}{{n + j}}} } \right)\\
&\quad+ 6{H_n}\left( {\sum\limits_{1 \le i < j \le r} {\frac{1}{{\left( {n + i} \right)\left( {n + j} \right)}}} } \right) + 6\sum\limits_{1 \le i < j < k \le r} {\frac{1}{{\left( {n + i} \right)\left( {n + j} \right)\left( {n + k} \right)}}}.\tag{2.24}
\end{align*}
Now, we use the notation ${W_{k,r}}\left( m \right)$ to stands for the sum
 \[{W_{k,r}}\left( m \right): =\left( {k - 1} \right)! \sum\limits_{n = 1}^\infty  {\frac{{\left[ {\begin{array}{*{20}{c}}
    {n + r + 1}  \\
    k  \\
 \end{array}} \right]}}{{\left( {n + r} \right)!{n^m}}}},\tag{2.25} \]
where $k\in \N,\ r\in \N_0$ and $m\in \mathbb{N} \setminus \{1\}$. By using the above notation, we obtain
\begin{align*}
&{W_{1,r}}\left( m \right) = \zeta \left( m \right), \\
&{W_{2,r}}\left( m \right) = \sum\limits_{n = 1}^\infty  {\frac{{{H_{n + r}}}}{{{n^m}}}} ,\\
&{W_{3,r}}\left( m \right) = \sum\limits_{n = 1}^\infty  {\frac{{H_{n + r}^2 - H_{n + r}^{\left( 2 \right)}}}{{{n^m}}}},\\
&{W_{4,r}}\left( m \right) = \sum\limits_{n = 1}^\infty  {\frac{{H_{n + r}^3 - 3{H_{n + r}}H_{n + r}^{\left( 2 \right)} + 2H_{n + r}^{\left( 3 \right)}}}{{{n^m}}}} .
\end{align*}
The following lemma will be useful in the development of the main theorem.
Noting that when $r=0$ and $k>1$ in (2.25), then using (1.11), which can be rewritten as
\begin{align*}
{W_{k,0}}\left( m \right):& = \left( {k - 1} \right)!\sum\limits_{n = 1}^\infty  {\frac{{\left[ {\begin{array}{*{20}{c}}
   {n + 1}  \\
   k  \\
\end{array}} \right]}}{{n!{n^m}}}} \\
& = \left( {k - 1} \right)!\left( {\sum\limits_{n = 1}^\infty  {\frac{{\left[ {\begin{array}{*{20}{c}}
   n  \\
   {k - 1}  \\
\end{array}} \right]}}{{n!{n^m}}} + \sum\limits_{n = 1}^\infty  {\frac{{\left[ {\begin{array}{*{20}{c}}
   n  \\
   k  \\
\end{array}} \right]}}{{n!{n^{m - 1}}}}} } } \right)\\
& = \left( {k - 1} \right)!\left( {\zeta \left( {m + 1,{{\left\{ 1 \right\}}^{k - 2}}} \right) + \zeta \left( {m,{{\left\{ 1 \right\}}^{k - 1}}} \right)} \right).\tag{2.26}
\end{align*}
On the other hand, the Aomoto-Drinfel¡¯d-Zagier formula reads
\[\sum\limits_{n,m = 1}^\infty  {\zeta \left( {m + 1,{{\left\{ 1 \right\}}^{n - 1}}} \right){x^m}{y^n} = 1 - \exp \left( {\sum\limits_{n = 2}^\infty  {\zeta \left( n \right)\frac{{{x^n} + {y^n} - {{\left( {x + y} \right)}^n}}}{n}} } \right)} ,\tag{2.27}\]
which implies that for any $m,\ n\in \N$, the multiple zeta value ${\zeta \left( {m + 1,{{\left\{ 1 \right\}}^{n - 1}}} \right)}$ can be
represented as a polynomial of zeta values with rational coefficients, and we have the duality formula
\[\zeta \left( {n + 1,{{\left\{ 1 \right\}}^{m - 1}}} \right) = \zeta \left( {m + 1,{{\left\{ 1 \right\}}^{n - 1}}} \right).\]
In particular, one can find explicit formulas for small weights $w:=n+m$.
\[\begin{array}{l}
 \zeta \left( {2,{{\left\{ 1 \right\}}^m}} \right) = \zeta \left( {m + 2} \right), \\
 \zeta \left( {3,{{\left\{ 1 \right\}}^m}} \right) = \frac{{m + 2}}{2}\zeta \left( {m + 3} \right) - \frac{1}{2}\sum\limits_{k = 1}^m {\zeta \left( {k + 1} \right)\zeta \left( {m + 2 - k} \right)} . \\
 \end{array}\]
Hence, we know that for $m,\ k\in \N$, the sums ${W_{k,0}}\left( m \right)$  can be expressed as a rational linear combination of zeta values. For example, from \cite{FS1998,CX2016}
\begin{align*}
&{W_{2,0}}\left( m \right) = \frac{1}{2}\left\{ {\left( {m + 2} \right)\zeta \left( {m + 1} \right) - \sum\limits_{i = 1}^{m - 2} {\zeta \left( {m - i} \right)\zeta \left( {i + 1} \right)} } \right\},\\
&{W_{3,0}}\left( m \right) = m{W_{2,0}}\left( m+1 \right) - \frac{{m\left( {m + 1} \right)}}{6}\zeta \left( {m + 2} \right) + \zeta \left( 2 \right)\zeta \left( m \right).
\end{align*}
\begin{lem}(\cite{X2016}) For integers $k\in \N$ and $p\in \mathbb{N} \setminus \{1\}$,  then the following identity holds:
\[\left( {p - 1} \right)!\sum\limits_{n = 1}^\infty  {\frac{{\left[ {\begin{array}{*{20}{c}}
   {n + 1}  \\
   p  \\
\end{array}} \right]}}{{n!n\left( {n + k} \right)}}}  = \frac{1}{k}\left\{ {\left( {p - 1} \right)!\zeta (p) + \frac{{{Y_p}\left( k \right)}}{p} - \frac{{{Y_{p - 1}}\left( k \right)}}{k}} \right\},\tag{2.28}\]
where ${Y_k}\left( n \right) = {Y_k}\left( {{H _n},1!{H^{(2)} _n},2!{H^{(3)}_n}, \cdots ,\left( {r - 1} \right)!{H^{(r)} _n}, \cdots } \right)$, ${Y_k}\left( {{x_1},{x_2}, \cdots } \right)$ stands for the complete exponential Bell polynomial defined by (see \cite{L1974})
\[\exp \left( {\sum\limits_{m \ge 1}^{} {{x_m}\frac{{{t^m}}}{{m!}}} } \right) = 1 + \sum\limits_{k \ge 1}^{} {{Y_k}\left( {{x_1},{x_2}, \cdots } \right)\frac{{{t^k}}}{{k!}}}.\tag{2.29}\]
\end{lem}
 From the definition of the complete exponential Bell polynomial, we have
$${Y_1}\left( n \right) = {H_n},{Y_2}\left( n \right) = H_n^2 + {H^{(2)} _n},{Y_3}\left( n \right) =  H_n^3+ 3{H_n}{H^{(2)} _n}+ 2{H^{(3)} _n},$$
\[{Y_4}\left( n \right) = H_n^4 + 8{H_n}{H^{(3)} _n} + 6H_n^2{H^{(2)} _n} + 3(H^{(2)} _n)^2 + 6{H^{(4)} _n},\]
\[{Y_5}\left( n \right) = H_n^5 + 10H_n^3{H^{(2)} _n} + 20H_n^2{H^{(3)}_n} + 15{H_n}({H^{(2)}_n})^2 + 30{H_n}{H^{(4)} _n}+ 20{H^{(2)} _n}{H^{(3)} _n} + 24{H^{(5)} _n}.\]
In fact, ${Y_k}\left( n \right)$ is a rational linear combination of products of harmonic numbers. Putting $p=2,3,4$ in (2.28), we obtain the corollary.
\begin{cor} For integer $k>0$, we have
\begin{align*}
&\sum\limits_{n = 1}^\infty  {\frac{{{H_n}}}{{n\left( {n + k} \right)}}}  = \frac{1}{k}\left( {\frac{1}{2}H_k^2 + \frac{1}{2}{H^{(2)} _k} + \zeta \left( 2 \right) - \frac{{{H_k}}}{k}} \right),\tag{2.30}\\
&\sum\limits_{n = 1}^\infty  {\frac{{H_n^2 - {H^{(2)} _n}}}{{n\left( {n + k} \right)}}}  = \frac{1}{k}\left\{ {2\zeta \left( 3 \right) + \frac{{H_k^3 + 3{H_k}{H^{(2)} _k} + 2{H^{(3)} _k}}}{3} - \frac{{H_k^2 + {H^{(2)} _k}}}{k}} \right\},\tag{2.31}\\
&\sum\limits_{n = 1}^\infty  {\frac{{H_n^3 - 3{H_n}{H^{(2)} _n} + 2{H^{(3)} _n}}}{{n\left( {n + k} \right)}}}  = \frac{1}{k}\left\{ \begin{array}{l}
 \frac{{H_k^4 + 8{H_k}{H^{(3)} _k} + 6H_k^2{H^{(2)} _k}+ 3(H^{(2)} _k)^2 + 6{H^{(4)} _k}}}{4} \\
  - \frac{{H_k^3 + 3{H_k}{H^{(2)} _k} + 2{H^{(3)} _k}}}{k} + 6\zeta \left( 4 \right) \\
 \end{array} \right\}.\tag{2.32}
\end{align*}
\end{cor}
Hence, combining (2.22)-(2.25), (2.30) and (2.31), we deduce the following identities
\begin{align*}
{W_{2,r}}\left( m \right)& = \sum\limits_{n = 1}^\infty  {\frac{{{H_n}}}{{{n^m}}}}  + \sum\limits_{j = 1}^r {\sum\limits_{n = 1}^\infty  {\frac{1}{{{n^m}\left( {n + j} \right)}}} } \\
& = {W_{2,0}}\left( m \right)+ \sum\limits_{l = 1}^{m - 1} {{{\left( { - 1} \right)}^{l - 1}}\zeta \left( {m + 1 - l} \right)H_r^{\left( l \right)}}  + {\left( { - 1} \right)^{m - 1}}\sum\limits_{j = 1}^r {\frac{{{H_j}}}{{{j^m}}}},\tag{2.33}\\
{W_{3,r}}\left( m \right)& = \sum\limits_{n = 1}^\infty  {\frac{{H_n^2 - H_n^{\left( 2 \right)}}}{{{n^m}}}}  + 2\sum\limits_{j = 1}^r {\sum\limits_{n = 1}^\infty  {\frac{{{H_n}}}{{{n^m}\left( {n + j} \right)}}} }  + 2\sum\limits_{1 \le i < j \le r} {\sum\limits_{n = 1}^\infty  {\frac{1}{{{n^m}\left( {n + i} \right)\left( {n + j} \right)}}} } \\
&= {W_{3,0}}\left( m \right) + 2\sum\limits_{i = 1}^{m - 1} {{{\left( { - 1} \right)}^{i - 1}}H_r^{\left( i \right)}{W_{2,0}}\left( {m + 1 - i} \right)} \\
&\quad + 2\sum\limits_{l = 1}^{m - 1} {{{\left( { - 1} \right)}^{l - 1}}\zeta \left( {m + 1 - l} \right)\sum\limits_{1 \le i < j \le r} {\frac{1}{{{i^l}\left( {j - i} \right)}}} } \\
&\quad - 2\sum\limits_{l = 1}^{m - 1} {{{\left( { - 1} \right)}^{l - 1}}\zeta \left( {m + 1 - l} \right)\sum\limits_{1 \le i < j \le r} {\frac{1}{{{j^l}\left( {j - i} \right)}}} } \\
&\quad + {\left( { - 1} \right)^{m - 1}}\left\{ \begin{array}{l}
 \sum\limits_{j = 1}^r {\frac{{H_j^2 + H_j^{\left( 2 \right)}}}{{{j^m}}}}  + 2\zeta \left( 2 \right)H_r^{\left( m \right)} - 2\sum\limits_{j = 1}^r {\frac{{{H_j}}}{{{j^{m + 1}}}}}  \\
  + 2\sum\limits_{1 \le i < j \le r} {\frac{{{H_i}}}{{{i^m}\left( {j - i} \right)}}}  - 2\sum\limits_{1 \le i < j \le r} {\frac{{{H_j}}}{{{j^m}\left( {j - i} \right)}}}  \\
 \end{array} \right\},\tag{2.34}\\
 {W_{4,r}}\left( m \right) &= {W_{4,0}}\left( m \right) + 3\sum\limits_{j = 1}^r {\sum\limits_{n = 1}^\infty  {\frac{{H_n^2 - H_n^{\left( 2 \right)}}}{{{n^m}\left( {n + j} \right)}}} }  + 6\sum\limits_{1 \le i < j \le r}\sum\limits_{n=1}^\infty {\frac{{{H_n}}}{{{n^m}\left( {n + i} \right)\left( {n + j} \right)}}} \\
 &\quad + 6\sum\limits_{1 \le i < j < k \le r}\sum\limits_{n=1}^\infty {\frac{1}{{{n^m}\left( {n + i} \right)\left( {n + j} \right)\left( {n + k} \right)}}} ,\\
 &  = {W_{4,0}}\left( m \right) + 3\sum\limits_{l = 1}^{m - 1} {{{\left( { - 1} \right)}^{l - 1}}H_r^{\left( l \right)}{W_{3,0}}\left( {m + 1 - l} \right)} \\
 &\quad + 6\sum\limits_{l = 1}^{m - 1} {{{\left( { - 1} \right)}^{l - 1}}{W_{2,0}}\left( {m + 1 - l} \right)\sum\limits_{1 \le i < j \le r} {\frac{1}{{j - i}}\left( {\frac{1}{{{i^l}}} - \frac{1}{{{j^l}}}} \right)} } \\
 &\quad + {\left( { - 1} \right)^{m - 1}}\left\{ \begin{array}{l}
 6\zeta \left( 3 \right)H_r^{\left( m \right)} - 3\sum\limits_{j = 1}^r {\frac{{H_j^2 + H_j^{\left( 2 \right)}}}{{{j^{m + 1}}}}}  \\
  + \sum\limits_{j = 1}^r {\frac{{H_j^3 + 3{H_j}H_j^{\left( 2 \right)} + 2H_j^{\left( 3 \right)}}}{{{j^m}}}}  \\
 \end{array} \right\}\\
 &\quad + 6{\left( { - 1} \right)^{m - 1}}\sum\limits_{1 \le i < j \le r} {\frac{1}{{j - i}}\left\{ \begin{array}{l}
 \frac{1}{{{i^m}}}\left( {\frac{{H_i^2 + H_i^{\left( 2 \right)}}}{2} + \zeta \left( 2 \right) - \frac{{{H_i}}}{i}} \right) \\
  - \frac{1}{{{j^m}}}\left( {\frac{{H_j^2 + H_j^{\left( 2 \right)}}}{2} + \zeta \left( 2 \right) - \frac{{{H_j}}}{j}} \right) \\
 \end{array} \right\}} \\
 &\quad+ 6\sum\limits_{l = 1}^{m - 1} {{{\left( { - 1} \right)}^{l - 1}}\zeta \left( {m + 1 - l} \right)} \sum\limits_{1 \le i < j < k \le r} {\left\{ \begin{array}{l}
 \frac{1}{{{i^l}\left( {j - i} \right)\left( {k - i} \right)}} + \frac{1}{{{j^l}\left( {i - j} \right)\left( {k - j} \right)}} \\
  + \frac{1}{{{k^l}\left( {j - k} \right)\left( {i - k} \right)}} \\
 \end{array} \right\}} \\
 &\quad + 6{\left( { - 1} \right)^{m - 1}}\sum\limits_{1 \le i < j < k \le r} {\left\{ \begin{array}{l}
 \frac{{{H_i}}}{{{i^m}\left( {j - i} \right)\left( {k - i} \right)}} + \frac{{{H_j}}}{{{j^m}\left( {i - j} \right)\left( {k - j} \right)}} \\
  + \frac{{{H_k}}}{{{k^m}\left( {j - k} \right)\left( {i - k} \right)}} \\
 \end{array} \right\}}.\tag{2.35}
\end{align*}
Therefore, the sums of harmonic numbers $W_{k,r}(m)$, for $k=1,2,3,4$ have been successfully represented in terms of zeta values and harmonic numbers. In fact, the other case of $W_{k,r}(m)$ can be evaluated in a similar fashion. Next, we shall present a closed form evaluation of the following sum:
\[S\left( {k,m;p} \right): = \sum\limits_{n = 1}^\infty  {\frac{{h_n^{\left( m \right)}\left( k \right)}}{{{n^p}}}},\ p\geq m+1,\ k=2,3. \]
By using the definitions of $S(k,m;p)$ and $W_{k,r}(m)$, then combining (2.20) and (2.21), we obtain the following theorem.
\begin{thm} For positive integers $r$ and $p\geq r+1$, then the following identities hold:
\begin{align*}
&S\left( {2,r;p} \right) = \frac{1}{{2!\left( {r - 1} \right)!}}\sum\limits_{k = 1}^r {\left[ {\begin{array}{*{20}{c}}
   r  \\
   k  \\
\end{array}} \right]} \left\{ \begin{array}{l}
 {W_{3,r - 1}}\left( {p + 1 - k} \right) - 2{H_{r - 1}}{W_{2,r - 1}}\left( {p + 1 - k} \right) \\
  + \left( {H_{r - 1}^2 + H_{r - 1}^{\left( 2 \right)}} \right)\zeta \left( {p + 1 - k} \right) \\
 \end{array} \right\},\tag{2.36}\\
 &S\left( {3,r;p} \right) = \frac{1}{{3!\left( {r - 1} \right)!}}\sum\limits_{k = 1}^r {\left[ {\begin{array}{*{20}{c}}
   r  \\
   k  \\
\end{array}} \right]} \left\{ \begin{array}{l}
 {W_{4,r - 1}}\left( {p + 1 - k} \right) - 3{H_{r - 1}}{W_{3,r - 1}}\left( {p + 1 - k} \right) \\
  + 3\left( {H_{r - 1}^2 + H_{r - 1}^{\left( 2 \right)}} \right){W_{2,r - 1}}\left( {p + 1 - k} \right) \\
  - \left( {H_{r - 1}^3 + 3{H_{r - 1}}H_{r - 1}^{\left( 2 \right)} + 2H_{r - 1}^{\left( 3 \right)}} \right)\zeta \left( {p + 1 - k} \right) \\
 \end{array} \right\}.\tag{2.37}
 \end{align*}
\end{thm}
From (2.33)-(2.37), we know that the sums $S(2,r;p)$ and $S(3,r;p)$ can be evaluated in terms of harmonic numbers and zeta values whenever $p\geq r+1$. A simple example is as follows:
\[S\left( {2,2;3} \right) = 4\zeta \left( 5 \right) - 2\zeta \left( 2 \right)\zeta \left( 3 \right) + \frac{5}{4}\zeta \left( 4 \right) - 2\zeta \left( 3 \right) + \zeta \left( 2 \right).\]
It may also be possible to represent the sums $S\left( {k,m;p} \right)$ for $4\leq k\in \N$ in closed form, this work is currently under investigation. It does appear however, that
there is a difficulty with the representation of $W_{k,r}(p)$ for $5\leq k\in \N$ in closed form.\\

 \end{document}